\documentclass[orivec,runningheads]{llncs}
\usepackage{epsf}
\usepackage[utf8]{inputenc}
\usepackage{amsmath}
\usepackage{graphicx}
\usepackage{soul}
\usepackage{algorithm}
\usepackage[noend]{algpseudocode}
\usepackage{rotating}
\usepackage{multirow}
\usepackage{subfigure}
\usepackage{amssymb}
\usepackage{url}
\usepackage{multicol} 
\usepackage{booktabs} 

\usepackage[colorinlistoftodos]{todonotes}

\newcommand{\keywords}[1]{\par\addvspace\baselineskip
\noindent\keywordname\enspace\ignorespaces#1}

\begin{document}

\title{A Statistical Comparison of Objetive Functions for the Vehicle Routing Problem with Route Balancing}

\titlerunning{A Statistical Comparison of Evaluation Functions for the Vehicle Routing Problem with Route Balancing}
\author{Jairo Lozano\inst{1} \and Luis C. Gonz\'alez-Gurrola\inst{1} \and Eduardo Rodriguez-Tello\inst{2} \and Philippe Lacomme\inst{3}}

\institute{Universidad Autónoma de Chihuahua, Facultad de Ingeniería\\
Circuito No. 1, Campus Universitario 2, 31125 Chihuahua, Chih., México.\\
          \email{{p245336, lcgonzalez}@uach.mx}
          \and
          CINVESTAV-Tamaulipas, Km. 5.5 Carretera Victoria-Soto La Marina, 87130 Victoria Tamps., México.\\
          \email{ertello@tamps.cinvestav.mx}
          \and 
          Laboratoire d'Informatique (LIMOS, UMR CNRS 6158), Campus des Cézeaux, 63177, Aubière Cedex, France\\
          \email{placomme@sp.isima.fr}
}

\maketitle


\begin{abstract}
The Vehicle Routing Problem with Route Balancing (VRPRB) is a bi-objective version of the original Vehicle Routing Problem (VRP) in which, besides minimizing the total distance traveled by the vehicles involved, the balance among route loads is also pursued. Different objective functions (OFs) to achieve balanced route configurations have been proposed in the literature, however to the best of the authors' knowledge there is still no consensus on which OF is the most suitable one for addressing, through metaheuristics, this challenging multi-objective optimization problem.

This paper inquires into the effectiveness of seven different OFs for the VRPRB.
Their influence on the performance of a basic single-solution-based evolutionary algorithm
is analyzed by comparing the quality of the Pareto-approximations produced for a set of well-known benchmark instances. The obtained results indicate that studying alternative evaluation schemes for the VRPRB represents a highly valuable direction for future research which merits more attention.


\end{abstract}



\section{Introduction}

The Vehicle Routing Problem (VRP) is a classical combinatorial optimization problem with great relevance both in theoretical research and in practical matter. Derived from the Traveling Salesman Problem (TSP), the VRP has been mainly tackled as a multi-objective problem, for two or three objectives of interest, therefore, a set of non-dominated solutions, the so-called \textit{Pareto-approximation} (see Sect.~\ref{multiobjectiveProblems}) is required. The most common objectives that have been addressed in literature are minimizing the total length of the routes as well as the number of vehicles used. Frequently, variants of VRP add more objectives on top of one of the two aforementioned. One flavor of the VRP that has received much less attention than other variants is the Vehicle Routing Problem with Route Balancing (VRPRB), which looks for the balance between routes as well as the total length minimization.  

A balanced VRP solution is one in which all routes are as similar as possible regarding a given indicator such as tour length or vehicle load.
This kind of solutions are required for example, when it comes to workload distribution among vehicle drivers or when trying to even out the travel time of passengers served by a transportation service.
For the later case, a problem that very well represents this requirement is the School Bus Routing Problem (SBRP) \cite{Bektas2007,Kim2012,Schittekat2013} that considers the pick-up and delivery of students to schools. Within this kind of problems, important aspects to consider are efficiency, effectiveness and equity, as firstly proposed by Savas \cite{Savas1978}. Nonetheless, when considering this scenario is clear that equity, defined as to minimize the difference in the time that students spend on the bus when being transported, has not  been fully researched. For example, Park and Kim \cite{Park2010} rightly stated: {\it ``Equity has been neglected in evaluating the performance of a school bus service as well as public services. However, its growing importance has been recognized. To improve the equity of a school bus service, balancing the loads and travel times between buses should be considered.''}. The importance of balance in real-world problems then, is attracting more attention, and the interest in the study of measures to achieve it is compelling. 
  
In order to capture the balance of a set of routes, different objective functions (OFs) have been defined in the VRPRB literature.
In \cite{Schwarze2013}, the authors identify six different types of OFs and give a brief explanation of each one to conclude with the selection of the OF \textit{Min-max} (minimization of the maximum indicator value) for their experiments. In \cite{Matl2016}, the authors identify a set of axiomatic properties that should be satisfied by an ideal equity OF and use them to evaluate six common measures of equity, pointing out their properties and the properties of the resulting Pareto-optimal fronts. For their study, the authors create a set of small VRP instances, solvable to optimality, and compute all the feasible tours. The authors conclude that none of the OFs satisfies every desirable axiom and thus no OF is better than the others in all relevant aspects. 

This article tries to motivate discussion and further research in the understanding of different OFs designed to capture balance. A statistical analysis is proposed to evaluate the quality of the Pareto-approximations obtained by using seven different OFs within a basic single-solution-based evolutionary algorithm, called \textit{(1+1) EA}, over the instances proposed in \cite{Matl2016}. We are interested in finding out what differences appear when applying specific OFs to the same data, using the same algorithm, and even the same random seed.

The rest of this paper is organized as follows. Sect.~\ref{multiobjectiveProblems} presents a basic overview of multi-objective optimization concepts and formally defines the VRPRB. Relevant literature is revisited in Sect.~\ref{literatureReview}. The seven OFs for the VRPRB analyzed in this study are defined in Sect.~\ref{materialsAndMethods}, where the components of the single-solution-based evolutionary algorithm used for assessing them are also described in detail. Sect.~\ref{experiments} presents the experiments performed for comparing the selected OFs as well as the results achieved. Finally, the conclusions of this study are stated in Sect.~\ref{conclusionsAndFutureWork}.


\section{Vehicle Routing Problem with Route Balancing}
\label{multiobjectiveProblems}
The VRP with balance requirements is {\it de facto} a multi-objective optimization problem. With multi-objective problems, in contrast to mono-objective problems, the objective is to provide a set of non-dominated solutions or \textit{Pareto-approximation}. Without loss of generality, the concept of dominance can be defined as follows. Given two solutions, where each solution $i$ has a vector $V_i=\{ob^i_1, ob^i_2, \ldots, ob^i_m\}$ of $m$ objectives, a solution $i$ dominates solution $j$, written $i \prec j$,  if $\forall q \in \{1, 2, \ldots, m\}$, $ob^i_q \le ob^j_q$ and $\exists$ $q$ such as $ob^i_q < ob^j_q$\footnote{For the sake of explanation this is a minimization problem, but it could be also generalized, previous change of relational operators, for maximization.}.  Given the NP-Hard nature of this kind of problems (the VRP and its variants), the set of non-dominated and \textit{optimal} solutions, i.e. the \textit{Pareto-front}, is unknown most of the times.




\subsection{VRP with Balance Requirements}

The Vehicle Routing Problem with Route Balancing (VRPRB) considers a set $V=\{0, 1, 2, \ldots, N\}$ of $N+1$ vertices, where $V_0$ represents the depot/school/work center, and the subset $V'=V$ \textbackslash $\{V_0\}$ represents the customers nodes to be attended. Each node $V_i$ has two attributes, a pair of geographic coordinates $(x_i,y_i)$ and a demand that needs to be collected $d_i$. For the depot $V_0$ its demand $d_0=0$. The pick-up service is provided by a set $K=\{2, 3, 4, \ldots, M\}$ of trucks/buses, where each truck/bus $K_j$ has capacity $C_j \ge d_i$, $\forall i \in \{1, 2, 3, \ldots, N\}$. For any pair of nodes $V_o$ and $V_s$, there is an associated travel cost $l_{o,s}$. The first objective is to find a set of feasible routes $T=\{1,2,3,\ldots,P\}$ with the minimum total distance, meaning that each route $T_w$ has an associated variable $L_w$ which represents the sum of all $l_{o,s}$, given that $o,s \in R_w$. The second objective which represents balance is to make as similar as possible any pair $L_w$, and $L_q$ $\forall w,q \in T$. In order to further motivate discussion, a quantitative measure for route balance is not provided in this section.




\section{Literature Review}
\label{literatureReview}
Most of the works reported in literature tackle the VRPRB as a bi-objective problem, where the first objective is to minimize  the total length of routes, and the second is the balance among the routes. 
For a first sub-classification, there are two variables that have been used for balancing purposes, i.e. load-related ones and the length of each route. When the load of routes is considered, the work of Lee and  Ueng \cite{Lee1999}  proposed an Integer Programming model with one of its objectives being the balance of the work-load of the drivers, then looking for fairness in the work schedules. In another work, Ribeiro and Louren\c{c}o \cite{Ribeiro2001} proposed a MILP model to balance volume transportation, thus trying to equally distribute compensations among drivers, since the earnings were associated to the load.  
When distances of the routes are used as a variable indicator for balancing, most of the works have used the $Max-min$ formulation (see Sect.~\ref{objectiveFunctions}). The first works in this direction were the ones of Josefowiez et al. \cite{Josefowiez2002,Josefowiez2006},  where they approached this problem proposing a modification of the NSGA II algorithm based on parallelization, afterwards, an implementation of a strategy based on a Pareto approach \cite{Josefowiez2007} was presented. Other work \cite{Josefowiez2009}, proposed a MOEA algorithm for this problem based on the islands models. Lacomme et. al. \cite{Lacomme2015} introduced a Multi-Start strategy based on Path-Relinking, where they compared their algorithm with the one proposed by Josefowiez et al. \cite{Josefowiez2009} obtaining better results in the Balance criterion. 
Recently, Schwarze and Vo\ss~\cite{Schwarze2013} presented a very complete survey of OFs used for balance and their indicator variables in the context of  Skill VRP, where load balancing and resource utilization is sought.



\section{Materials and Methods}
\label{materialsAndMethods}
This section introduces the OFs under study, the algorithm and the two mutation operators implemented as well as a TSP-to-VRP transformation involved in the mutation process. A description of the heuristic used for initial solution creation and the metrics used for comparison purposes are also given.

\subsection{Objective Functions Under Study}
\label{objectiveFunctions}

Seven OFs for route balancing were chosen from the literature including some of the most widely used and some others that are not that common but have been taken into account in some of the latest studies. Among the most commonly used OFs we find \textit{Max-min} (also called \textit{range}), perhaps the most widely used, which tries to minimize the difference between the maximum and minimum indicator value; \textit{Min-max} which tries to minimize the maximum indicator value; and \textit{Var} which tries to minimize the variance of indicator values (some studies also use the standard deviation for this matter). In \cite{Schwarze2013}, the authors also identify an objective function trying to minimize the relative deviation of indicator values from maximum value, namely \textit{Rel} and another trying to minimize the cumulative difference (over all routes) between the indicator value for the given route and the minimum indicator value, namely \textit{All-min}. In a more recent study Matl et. al. \cite{Matl2016} analyzed two more OFs: \textit{MAD} (Mean Absolute Deviation), defined as the mean absolute difference between each indicator value and the mean of indicator values; and the \textit{Gini Coefficient} or \textit{Gini}, one of the most widely used measure for inequality studies in general.

A formal definition of each OF under study is presented:
\newline
\newline
\noindent\textbf{All-min}\newline
\begin{equation}
	min \sum_{t\in T}(l_{t} - min_{u\in T}\mathit{l}_{u})
\end{equation}
\textbf{Max-min}\newline
\begin{equation}
	min (max_{u\in T}\mathit{l}_{u} - min_{u\in T}\mathit{l}_{u})
\end{equation}
\textbf{Min-max}\newline
\begin{equation}
	min \ max_{u\in T}\mathit{l}_{u}
\end{equation}
\textbf{Rel}\newline
\begin{equation}
	min \ \frac{1}{\left | T \right |} \ \sum_{t\in T} \left(\frac{max_{u\in T}l_{u} - l_{t}}{max_{u\in T} l_{u}}\right)
\end{equation}
\textbf{Var}\newline
\begin{equation}
min \ \left( \frac{\sum_{t\in T} l_{t}^{2}}{\left | T \right |} \ - \left( \frac{\sum_{t\in T} l_{t}}{\left | T \right |}\right)^{2} \right)
\end{equation}
\textbf{MAD}\newline
\begin{equation}
min \ \frac{1}{\left | T \right |} \ \sum_{t\in T} \left | l_{i} - \bar{l} \right |
\end{equation}
\textbf{Gini}\newline
\begin{equation}
min \ \frac{1}{2n^{2}\bar{l}} \ \sum_{t\in T} \ \sum_{{t}'\in T} \left | l_{t} - l_{{t}'} \right |
\end{equation}

Where $T$ stands for the set of routes of a given solution and $l$, the indicator value, was chosen to be the tour length. 


\subsection{(1+1) Evolutionary Algorithm with External Archive}
\label{sectionAlgorithm}

In order to find out if using different OFs for balance leads to solutions of different quality when used within a meta-heuristic, the simplest version of an evolutionary algorithm, the (1+1) EA, was implemented, which consists in the iterative mutation of a single individual $sol$ created initially at random. The individual $sol$ is only replaced by the mutated one $sol'$ when the latter is at least as good as the former. But, considering the bi-objective nature of the VRPRB,
a variant as presented by \cite{Garza-Fabre2012} was chosen, which includes an external archive to store the non-dominated solutions found during the evolution process in such a way that the mutated individual $sol'$ is only accepted if it is non-dominated by any individual already in the archive. When a new individual $sol'$ is accepted, all individuals dominated by $sol'$, as well as those mapping to the same objective vector, are removed from the non-dominated archive.
The (1+1)EA with external archive is formally described in Algorithm~\ref{alg:1+1ea}.

\begin{algorithm}
\caption{(1+1) EA with external Archive}\label{alg:1+1ea}
\begin{algorithmic}[1]
\Procedure{solveInstance}{}
\State $sol\gets$  randomly generated solution
\State $A \gets \left\{sol\right\}$
\State $iter \gets 0$
\While{$iter < MAX\_ITER$}
\State $sol' \gets $ mutationOperator($sol$)
\If {$\nexists\hat{sol}\in A: \hat{sol} \prec sol'$}
\State $A \gets \left\{ \hat{sol} \in A:sol' \nprec \hat{sol}  \wedge f(\hat{sol}) \neq f(sol') \right\} \cup \left\{ sol' \right\} $
\State $sol \gets sol'$
\EndIf
\State $iter \gets iter + 1$
\EndWhile\label{euclidendwhile}
\EndProcedure
\end{algorithmic}
\end{algorithm}

The selection of such a simple algorithm was in order to avoid any possible bias towards the benefit of any specific objective function. $MAX\_ITER$ was fixed to 30,000, which is a large enough number of iterations for the algorithm to converge.
Also in order to get impartial results, two different mutation operators were developed: \textit{Swap} mutation and \textit{Reverse} mutation. Both mutation operators make use of a search space transformation from VRP to TSP and vice versa before and after the solution disturbance.
In the next paragraphs the heuristic used to generate the initial solution, the mutation operators and the TSP-VRP transformation process are described.


\subsubsection{Heuristic for Constructing the Initial Solution}
For creating the initial solution, a client permutation list (containing all clients) is randomly generated.
Then, clients are inserted one by one into the first route until no more clients can be included due to the vehicle capacity constraint. When no more clients can be included into the first route, a second route is created to insert the following clients. The process stops until all clients have been inserted into a trip.

\subsubsection{VRP-to-TSP Transformation}
Because (1+1) EA uses only a mutation operator it is difficult to get an appropriate balance between exploration (searching on a large portion of the search space trying to find more promising areas) and exploitation (searching on a limited and promising region of the search space trying to improve a high quality solution). In order to achieve this balance between exploration and exploitation, a VRP-to-TSP transformation, before the solution perturbation, is used. This is achieved by concatenating all the individual routes into a TSP permutation. Operating over the TSP version of the solution allows the mutation operators to work inside a single route as well as between different routes, thus favoring both exploitation and exploration, respectively. Once the perturbation has been done, the TSP solution is turned back into a VRP solution in a similar way to that used to construct the initial solution.


\subsubsection{Swap Mutation Operator}
\label{swapMutation}
First, the VRP solution is transformed into a TSP permutation. Then, a random number $n$ of clients to swap is chosen (for this study, $n$ was bounded between 2 and 4). Then, $n$ positions from the TSP permutation are chosen and swapped randomly. The resulting TSP permutation is finally transformed back into a VRP solution.


\subsubsection{Reverse Mutation}
\label{reverseMutation}
The original VRP solution is encoded as a TSP solution and a randomly selected sub-path of the TSP is reversed. Finally, the resulting TSP permutation is transformed back into a VRP solution.

\bigskip
Note that neither the (1+1) EA, nor the mutation operators (\textit{Swap} and \textit{Reverse}) include any strategy to make individual routes optimal in length (optimal TSPs), thus solutions could be artificially balanced by increasing the length of small routes within a solution, as stressed by Jozefowiez et. al. in \cite{Josefowiez2009}. 



\subsection{Hypervolume to Measure the Pareto-approximation Quality}

When applied to a multi-objective problem, a meta-heuristic would generally lead to a set of non-dominated solutions known as Pareto-approximation. Several metrics could be used to measure the quality of a \textit{Pareto-approximation}. One of the most popular metrics is the hypervolume metric H(A,{\bf z}) \cite{zitzler1998}, which captures the area dominated by an approximation set $A$ and a reference point {\bf z}. The larger the hypervolume, the better the convergence of a given solution set to the optimum set (Pareto-front).
 
A more formal definition for the hypervolume is as follows. Considering a two objective dimensional space $f(x)=(f_1(x), f_2(x))$, each solution $x_i \in A$ forms a rectangle defined by its coordinates $(f_1(x_i),f_2(x_i))$ and the coordinates of the reference point $(z_1,z_2)$. The hypervolume is calculated as the union of all the rectangles defined between all the points in $A$ and {\bf z}. For minimization problems, one typical decision is to choose the worst possible value on each objective as the coordinates of the reference point. For this study, the reference point is computed as follows. For the balance, a solution with only two routes is created, the first one including all the clients except the last one, which is included into the second route. For the distance, a solution with one route per client is created. This way a solution with extremely poor objective values is obtained. It is important to point out that, in order to get a hypervolume measure equally influenced by both objectives, the objective values of each solution in the Pareto-approximation are normalized to a 0-1 range.





\section{Experiments}
\label{experiments}
In this section, we describe the experiments that were performed in order to statistically compare the OFs presented in Sect.~\ref{objectiveFunctions}. First, the set of VRP instances recently proposed in \cite{Matl2016} was selected. The set is composed of 60 small VRP instances, solvable to optimality, which are created as follows: based on the instances 1-10 of Christofides et. al. \cite{Christofides1979}, the authors take sets of $n = 14$ customers to create a subset of instances by varying the number of vehicles $\left | T \right |$ between 2 and 5, and the vehicle capacity $q$  to $\left[n/ \left | T \right | \right]$ or $\left[n/ \left | T \right |  \right] + 1$. All the customers demands are set to 1.

The algorithm described in Sect.~\ref{sectionAlgorithm} was coded in JAVA and run on a PC with 2.50 GHz CPU and 6 GB of RAM. Two variants are considered: the first one using the \textit{Swap} mutation operator and the other using the \textit{Reverse} mutation operator. 
Both versions are benchmarked on the 60 VRP instances using each of the seven OFs under study, so 14 different configurations (\textit{mutation operator-OF}) are tested.
Each variant is run 30 times per instance and each run goes for 30,000 iterations. At the end of each run, an archive with a set of non-dominated solutions is obtained, and its hypervolume is computed.


\begin{table}[tbh]
\centering
\caption{Cumulative scores obtained by each compared OF over 60 instances using two different mutation operators.} \label{tab:Results} 
\setlength{\tabcolsep}{3 pt}
\renewcommand{\arraystretch}{0.95}
\begin{scriptsize}
\begin{tabular}{*{16}{r}}
	\hline
	\noalign{\smallskip}
	\multicolumn{1}{c}{} & \multicolumn{7}{c}{Reverse mutation} & & \multicolumn{7}{c}{Swap mutation} \\
	\cmidrule{2-8} \cmidrule{10-16}
	\multicolumn{1}{c}{\rotatebox{65}{Instance}} & \multicolumn{1}{c}{\rotatebox{65}{All-min}} & \multicolumn{1}{c}{\rotatebox{65}{Max-min}} & \multicolumn{1}{c}{\rotatebox{65}{Min-max}} & \multicolumn{1}{c}{\rotatebox{65}{Rel}} & \multicolumn{1}{c}{\rotatebox{65}{Var}} & \multicolumn{1}{c}{\rotatebox{65}{MAD}} & \multicolumn{1}{c}{\rotatebox{65}{Gini}} & & \multicolumn{1}{c}{\rotatebox{65}{All-min}} & \multicolumn{1}{c}{\rotatebox{65}{Max-min}} & \multicolumn{1}{c}{\rotatebox{65}{Min-max}} & \multicolumn{1}{c}{\rotatebox{65}{Rel}} & \multicolumn{1}{c}{\rotatebox{65}{Var}} & \multicolumn{1}{c}{\rotatebox{65}{MAD}} & \multicolumn{1}{c}{\rotatebox{65}{Gini}} \\
	\hline
	\noalign{\smallskip}
1 & 1 & 1 & -6 & 0 & \textbf{2} & 1 & 1 && \textbf{1} & \textbf{1} & -6 & \textbf{1} & \textbf{1} & \textbf{1} & \textbf{1} \\
2 & -2 & -2 & -6 & -2 & \textbf{4} & \textbf{4} & \textbf{4} && -2 & -2 & -6 & -2 & \textbf{4} & \textbf{4} & \textbf{4} \\
3 & -2 & -2 & -6 & -2 & 4 & \textbf{5} & 3 && -2 & -2 & -6 & -2 & \textbf{4} & \textbf{4} & \textbf{4} \\
4 & -1 & -1 & -6 & -1 & -1 & \textbf{5} & \textbf{5} && -1 & -2 & -6 & -3 & \textbf{4} & \textbf{4} & \textbf{4} \\
5 & -3 & -1 & -6 & -2 & 4 & 3 & \textbf{5} && -2 & -2 & -6 & -2 & \textbf{4} & \textbf{4} & \textbf{4} \\
6 & -4 & -1 & -6 & -1 & \textbf{4} & \textbf{4} & \textbf{4} && -3 & -2 & -6 & -1 & \textbf{4} & \textbf{4} & \textbf{4} \\
7 & 1 & 0 & -6 & \textbf{2} & 1 & 1 & 1 && 1 & 0 & -6 & 1 & 1 & 1 & \textbf{2} \\
8 & -2 & -1 & -6 & -2 & 3 & \textbf{6} & 2 && -2 & -2 & -6 & -2 & \textbf{4} & \textbf{4} & \textbf{4} \\
9 & -2 & -2 & -6 & -2 & \textbf{4} & \textbf{4} & \textbf{4} && -2 & -2 & -6 & -2 & \textbf{4} & \textbf{4} & \textbf{4} \\
10 & -3 & -2 & -5 & -2 & \textbf{4} & \textbf{4} & \textbf{4} && -2 & -2 & -6 & -2 & 3 & 4 & \textbf{5} \\
11 & -2 & -2 & -6 & -2 & \textbf{4} & \textbf{4} & \textbf{4} && -2 & -2 & -6 & -2 & \textbf{4} & \textbf{4} & \textbf{4} \\
12 & -2 & -1 & -6 & -3 & \textbf{4} & \textbf{4} & \textbf{4} && -2 & -2 & -6 & -2 & \textbf{5} & 2 & \textbf{5} \\
13 & 0 & 0 & -1 & 0 & -1 & 0 & \textbf{2} && 0 & \textbf{1} & -4 & \textbf{1} & 0 & \textbf{1} & \textbf{1} \\
14 & -2 & -2 & -6 & -2 & \textbf{4} & \textbf{4} & \textbf{4} && -2 & -2 & -6 & -2 & 4 & 3 & \textbf{5} \\
15 & -3 & -3 & -3 & -3 & 3 & 4 & \textbf{5} && -2 & -2 & -6 & -2 & \textbf{4} & \textbf{4} & \textbf{4} \\
16 & -2 & -2 & -5 & -3 & \textbf{4} & \textbf{4} & \textbf{4} && -2 & -2 & -6 & -2 & \textbf{4} & \textbf{4} & \textbf{4} \\
17 & -2 & -2 & -6 & -2 & \textbf{4} & \textbf{4} & \textbf{4} && -2 & -1 & -6 & -3 & \textbf{4} & \textbf{4} & \textbf{4} \\
18 & -2 & -2 & -6 & -2 & \textbf{4} & \textbf{4} & \textbf{4} && -2 & -2 & -6 & -2 & \textbf{4} & \textbf{4} & \textbf{4} \\
19 & 0 & 2 & -6 & -2 & 2 & \textbf{3} & 1 && 1 & \textbf{2} & -6 & 0 & 1 & 1 & 1 \\
20 & -2 & -2 & -6 & -2 & \textbf{4} & \textbf{4} & \textbf{4} && -2 & -2 & -6 & -2 & \textbf{4} & \textbf{4} & \textbf{4} \\
21 & -2 & -2 & -6 & -2 & \textbf{4} & \textbf{4} & \textbf{4} && -1 & -1 & -6 & 0 & 1 & \textbf{4} & 3 \\
22 & -2 & -2 & -6 & -2 & 4 & 3 & \textbf{5} && -2 & -2 & -6 & -2 & \textbf{4} & \textbf{4} & \textbf{4} \\
23 & -4 & -3 & -4 & -1 & \textbf{4} & \textbf{4} & \textbf{4} && -2 & -3 & -6 & -1 & \textbf{4} & \textbf{4} & \textbf{4} \\
24 & -3 & -2 & -3 & -4 & \textbf{4} & \textbf{4} & \textbf{4} && -3 & 0 & -6 & -3 & \textbf{5} & 4 & 3 \\
25 & -1 & -1 & -5 & -2 & \textbf{5} & -1 & \textbf{5} && 0 & 1 & -6 & \textbf{3} & 1 & 0 & 1 \\
26 & -2 & -2 & -6 & -2 & 3 & \textbf{5} & 4 && -2 & -2 & -6 & -2 & \textbf{4} & \textbf{4} & \textbf{4} \\
27 & -2 & -2 & -6 & -2 & \textbf{4} & \textbf{4} & \textbf{4} && -2 & -2 & -6 & -2 & \textbf{4} & \textbf{4} & \textbf{4} \\
28 & -2 & -2 & -6 & -2 & \textbf{4} & \textbf{4} & \textbf{4} && -2 & -2 & -6 & -2 & \textbf{4} & \textbf{4} & \textbf{4} \\
29 & -2 & -2 & -6 & -2 & 3 & \textbf{5} & 4 && -2 & -2 & -6 & -2 & \textbf{4} & \textbf{4} & \textbf{4} \\
30 & -2 & -2 & -6 & -2 & \textbf{5} & 3 & 4 && -2 & -1 & -6 & -3 & \textbf{4} & \textbf{4} & \textbf{4} \\
31 & 0 & 1 & -6 & 1 & 1 & \textbf{2} & 1 && \textbf{1} & \textbf{1} & -6 & \textbf{1} & \textbf{1} & \textbf{1} & \textbf{1} \\
32 & -2 & -2 & -6 & -2 & \textbf{4} & \textbf{4} & \textbf{4} && -2 & -2 & -6 & -2 & \textbf{4} & \textbf{4} & \textbf{4} \\
33 & -2 & -2 & -6 & -2 & \textbf{4} & \textbf{4} & \textbf{4} && -2 & -2 & -6 & -2 & \textbf{4} & \textbf{4} & \textbf{4} \\
34 & -2 & -2 & -6 & -2 & \textbf{4} & \textbf{4} & \textbf{4} && -2 & -2 & -6 & -2 & \textbf{4} & \textbf{4} & \textbf{4} \\
35 & -2 & -2 & -6 & -2 & 4 & 3 & \textbf{5} && -2 & -2 & -6 & -2 & \textbf{4} & \textbf{4} & \textbf{4} \\
36 & -1 & -2 & -6 & -3 & 4 & 3 & \textbf{5} && -1 & -2 & -6 & -3 & \textbf{4} & \textbf{4} & \textbf{4} \\
37 & \textbf{1} & \textbf{1} & -6 & \textbf{1} & \textbf{1} & \textbf{1} & \textbf{1} && \textbf{2} & 1 & -6 & 1 & 1 & 0 & 1 \\
38 & -3 & -3 & -3 & -3 & \textbf{4} & \textbf{4} & \textbf{4} && -2 & -2 & -6 & -2 & \textbf{4} & \textbf{4} & \textbf{4} \\
39 & -1 & -2 & -6 & -3 & \textbf{4} & \textbf{4} & \textbf{4} && -1 & -2 & -6 & -1 & 2 & \textbf{4} & \textbf{4} \\
40 & -1 & -1 & -6 & -1 & -1 & \textbf{5} & \textbf{5} && -2 & -2 & -6 & -2 & \textbf{4} & \textbf{4} & \textbf{4} \\
41 & -2 & -2 & -5 & -3 & \textbf{4} & \textbf{4} & \textbf{4} && -2 & -2 & -6 & -2 & \textbf{4} & \textbf{4} & \textbf{4} \\
42 & -1 & -2 & -6 & -3 & \textbf{4} & \textbf{4} & \textbf{4} && -1 & -2 & -6 & -3 & \textbf{4} & \textbf{4} & \textbf{4} \\
43 & \textbf{1} & \textbf{1} & -6 & \textbf{1} & \textbf{1} & \textbf{1} & \textbf{1} && 1 & 1 & -6 & 1 & \textbf{2} & 1 & 0 \\
44 & -4 & -3 & -2 & -3 & \textbf{4} & \textbf{4} & \textbf{4} && -3 & -3 & -3 & -3 & \textbf{4} & \textbf{4} & \textbf{4} \\
45 & -2 & -2 & -6 & -2 & \textbf{4} & \textbf{4} & \textbf{4} && -2 & -2 & -6 & -2 & \textbf{4} & \textbf{4} & \textbf{4} \\
46 & -3 & -3 & -3 & -3 & \textbf{4} & \textbf{4} & \textbf{4} && -2 & -2 & -6 & -2 & \textbf{4} & \textbf{4} & \textbf{4} \\
47 & -3 & -3 & -3 & -3 & \textbf{4} & \textbf{4} & \textbf{4} && -2 & -2 & -6 & -2 & 3 & 4 & \textbf{5} \\
48 & -2 & -2 & -6 & -2 & \textbf{4} & \textbf{4} & \textbf{4} && -2 & -2 & -6 & -2 & \textbf{4} & \textbf{4} & \textbf{4} \\
49 & \textbf{1} & \textbf{1} & -6 & \textbf{1} & \textbf{1} & \textbf{1} & \textbf{1} && \textbf{1} & \textbf{1} & -6 & \textbf{1} & \textbf{1} & \textbf{1} & \textbf{1} \\
50 & -2 & -2 & -6 & -2 & \textbf{4} & \textbf{4} & \textbf{4} && -2 & -2 & -6 & -2 & \textbf{4} & \textbf{4} & \textbf{4} \\
51 & -2 & -2 & -6 & -2 & \textbf{4} & \textbf{4} & \textbf{4} && -2 & -2 & -6 & -2 & \textbf{4} & \textbf{4} & \textbf{4} \\
52 & 0 & 0 & -6 & 0 & -4 & 4 & \textbf{6} && -2 & -3 & -6 & -1 & \textbf{4} & \textbf{4} & \textbf{4} \\
53 & -2 & -2 & -6 & -2 & \textbf{4} & \textbf{4} & \textbf{4} && -2 & -2 & -6 & -2 & \textbf{4} & \textbf{4} & \textbf{4} \\
54 & -1 & -2 & -6 & -3 & \textbf{5} & 2 & \textbf{5} && -2 & -2 & -6 & -2 & \textbf{4} & \textbf{4} & \textbf{4} \\
55 & 0 & -3 & -6 & 2 & 2 & \textbf{3} & 2 && \textbf{1} & \textbf{1} & -6 & \textbf{1} & \textbf{1} & \textbf{1} & \textbf{1} \\
56 & -2 & -2 & -6 & -2 & \textbf{4} & \textbf{4} & \textbf{4} && -2 & -2 & -6 & -2 & \textbf{4} & \textbf{4} & \textbf{4} \\
57 & -2 & -2 & -6 & -2 & \textbf{4} & \textbf{4} & \textbf{4} && -2 & -2 & -6 & -2 & 3 & \textbf{6} & 3 \\
58 & -2 & -2 & -6 & -2 & \textbf{4} & \textbf{4} & \textbf{4} && -2 & -2 & -6 & -2 &\textbf{4} & \textbf{4} & \textbf{4} \\
59 & -3 & -1 & -6 & -2 & \textbf{4} & \textbf{4} & \textbf{4} && -2 & -2 & -6 & -2 & \textbf{4} & \textbf{4} & \textbf{4} \\
60 & -2 & -1 & -6 & -3 & 3 & \textbf{6} & 3 && -3 & 0 & -6 & -3 & \textbf{4} & \textbf{4} & \textbf{4} \\
	\hline
\end{tabular} 
\end{scriptsize}
\end{table}

Thirty Pareto-approximations were obtained per instance for each of the 14 algorithm variants (mutation and objective function combinations). Since every objective function works on a different objective space for balancing purposes, in order to statistically compare the quality of the obtained sets, a common space where all the objective functions can be fairly compared should be considered. Thus, the $Max-min$ objective space was selected since it is the most widely used in the literature. So each solution in the Pareto-approximations is re-evaluated using the $Max-Min$ OF, then the hypervolume of each Pareto-approximation over a unified objective space can be computed, allowing a fair comparison.

To statistically compare the vectors composed of 30 hypervolumes obtained by using each of the seven OFs, a $t$-test (two-sample, two-tailed, unequal variance) is applied between pairs of vectors. The results presented in this section are based on a scoring scheme that works as follows. When comparing two hypervolume vectors, corresponding to two different OFs, if the $t$-test shows that differences between both vector averages are significant at a 95\% confidence level, we add +1 to the score of the OF with larger average and -1 in detriment to the other OF score. If there is no statistically significant difference between the two vector averages, 0 is added for both OF scores.
    
In Table~\ref{tab:Results} the scores obtained by each OF over 60 instances when using two different mutation operators are presented. The highest scores per instance are marked in bold for each mutation operator. For example, for instance number 35 and the \textit{Reverse} mutation, \textit{MAD} gets an accumulated score of 3, given that it overcomes OFs \textit{All-min}, \textit{Max-min}, \textit{Min-max} and \textit{Rel} (in terms of the $t$-test analysis previously described), so +4 is added to its score, but it is also overcame by \textit{Gini} and indifferent with \textit{Var}, so -1 and 0 is added to obtain the final score of 3.

A clear disadvantage of the \textit{Min-max} OF is noticed, since only negative scores for both mutation operators are obtained. It would be of great interest to investigate if this also happens when using a more powerful meta-heuristic.

It is also important to note that for some instances, several OFs get the highest score, meaning that no OF dominates all of the others; e.g. instances 10, 16 or 46 for the \textit{Reverse} mutation and instances 4, 8 or 60 for the \textit{Swap} mutation, where \textit{Var}, \textit{MAD} and \textit{Gini} got the highest scores. And there are also some instances in which no OF appears to have a significant advantage as is the case for instances 37 and 43 when using the \textit{Reverse} mutation and instances 1 and 49 when using the \textit{Swap} mutation.

It is interesting how the scores are very similar in both tables even when they were obtained by using quite different mutation operators.

\begin{figure}[tb]
	\includegraphics[width=123mm]{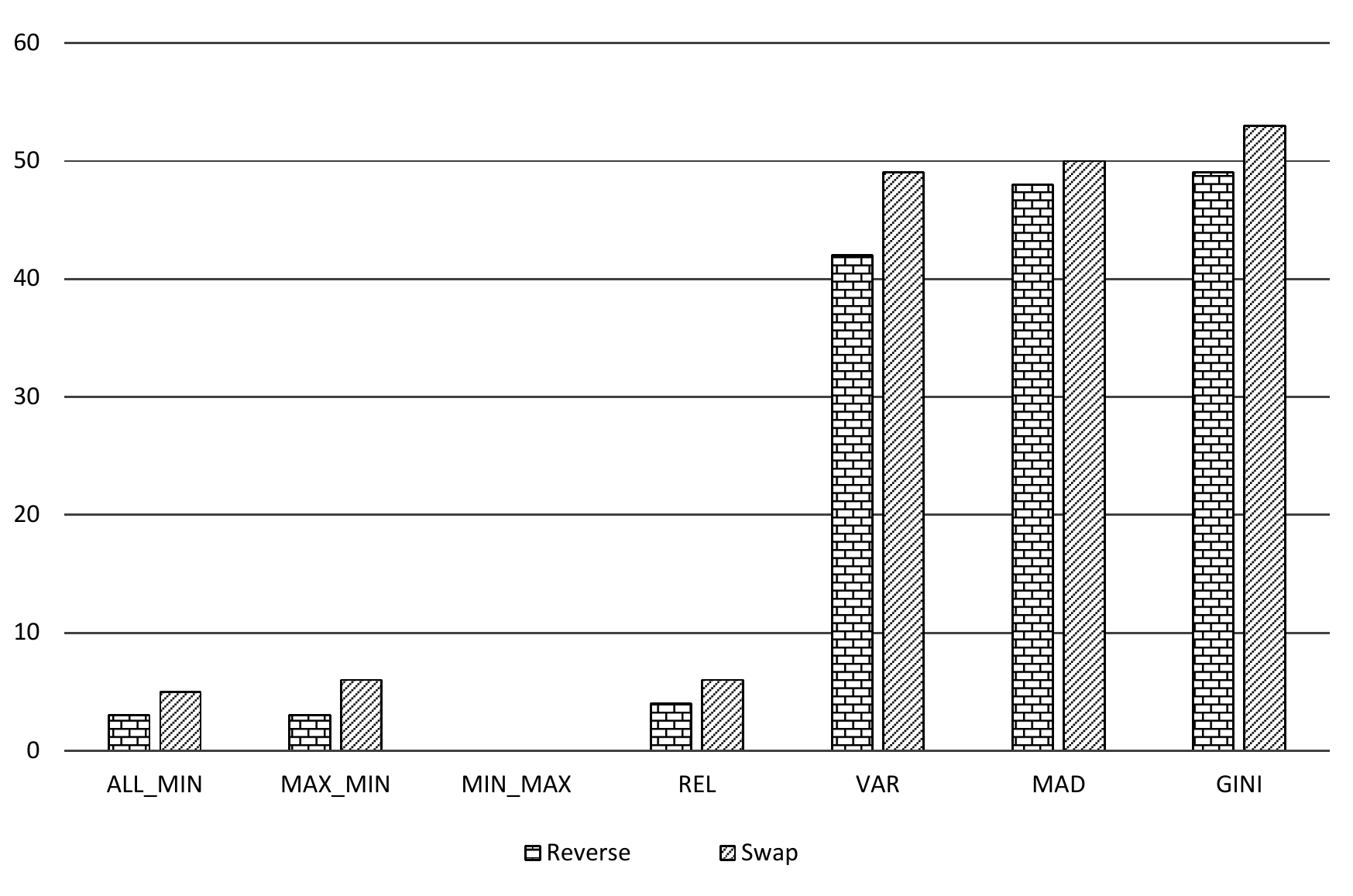}
    \caption{Number of times a given objective function got the highest score over the 60 instances using the two mutation operators.}\label{fig:histogram}
\end{figure}

In Fig.~\ref{fig:histogram} a histogram is presented showing how many times a given OF got the highest score after the $t$-test evaluations over the 60 instances using the two different mutation operators.

Here some facts attract attention. A clear domination of OFs \textit{Var}, $MAD$ and $Gini$ is observed, meaning that these three OFs were the ones with statistically higher hypervolume vector-averaged over the 60 instances. It is noticeable that \textit{Max-min}, being perhaps the most widely used OF for balance in the VRPRB literature, has been overcame most of the times by these three OF. The same goes for \textit{Min-max} which was unable to beat the rest of the OFs even once. 

It is also interesting that, even when the two mutation operators work quite differently, the scores obtained by the different OFs are very similar, for this matter it would also be of interest to investigate if the scores behave similar when using a more powerful meta-heuristic.

Using two quite different and free of bias algorithms (mutation operators) under exactly the same parameters over a set of 60 different VRP instances, and observing very similar behaviors when using the different OFs under study, pushes us into considering the next. Rather than being correlated to the OF alone, the quality of the Pareto-approximation obtained when using different OFs for balance must be somehow correlated with some properties of the VRP instance itself, like the distribution of the clients or the number of vehicles. This constitutes a very interesting topic for further research.
	



\section{Conclusions and Future Work}
\label{conclusionsAndFutureWork}
A statistical comparison was made between quality (in terms of hypervolume) of the Pareto-approximations obtained when using seven different objective functions for balance in a (1+1) EA with two different mutation operators for VRPRB over a set of a state-of-the-art instances proposed by \cite{Matl2016}. The size of the instances as well as the type of algorithm and the mutation operators were chosen in order to avoid any bias towards the benefit of any particular objective function.

Experiments clearly demonstrate three of the objective functions under study, namely \textit{Var}, \textit{MAD} and \textit{Gini}, leading to Pareto-approximations of higher quality than the rest of the objective functions, thus overcoming two of the most widely used in the VRPRB literature: \textit{Max-min} and \textit{Min-max}, and two other also considered: \textit{All-min} and \textit{Rel}.

Moreover, similar behaviors were observed when using the two different mutation operators: a clear advantage of \textit{Var}, \textit{MAD} and \textit{Gini} objective functions when it comes to how many times the objective function got the highest score in the $t$-test performed by instance; poor quality Pareto-approximations for the rest of the objective functions, highlighting \textit{Max-min} and \textit{Min-max} which are two of the most widely used objective functions in the VRPRB literature.

This pushes us into considering that the quality of the Pareto-approximations when using different objective functions does not depend on the algorithm used nor in the objective function alone, but in a possible correlation between some of the VRP instance properties and the objective function used. Further research on this matter is of great interest.

In a future study it would be interesting to include a more powerful meta-heuristic and tackle some more complex VRP instances to determine if behaviors among algorithms when using different objective functions for balance are still similar and thus reinforcing the hypothesis stated in the previous paragraph.





\bibliographystyle{ieeetr}
\bibliography{paper}
\end{document}